\documentclass[11pt]{article}

\usepackage{fullpage}
\usepackage{xypic}
\usepackage{amsgen}
\usepackage{amsmath}
\usepackage{amstext}
\usepackage{amsbsy}
\usepackage{amsopn}
\usepackage{amsfonts}
\usepackage{amssymb}
\usepackage{eepic}
\usepackage{graphicx}
\usepackage{epsf}
\usepackage{pstricks}
\xyoption{all}

\def\Box{\square}

\def\mapright#1{\smash{\mathop{\longrightarrow}\limits^{#1}}}
\def\tra#1{\smash{\mathop{\mid\kern
-1pt\joinrel\relbar\joinrel\relbar}\limits^{*}_{#1}}}
\def\longtra#1{\smash{\mathop{\mid\kern
-1pt\joinrel\relbar\joinrel\relbar\joinrel\relbar}\limits^{*}_{#1}}}
\def\vlongtra#1{\smash{\mathop{\mid\kern
-1pt\joinrel\relbar\joinrel\relbar\joinrel\relbar\joinrel\relbar}\limits^{*}_{#1}}}
\def\vvlongtra#1{\smash{\mathop{\mid\kern
-1pt\joinrel\relbar\joinrel\relbar\joinrel\relbar\joinrel\relbar\joinrel\relbar}\limits^{*}_{#1}}}
\def\vvvlongtra#1{\smash{\mathop{\mid\kern
-1pt\joinrel\relbar\joinrel\relbar\joinrel\relbar\joinrel\relbar\joinrel\relbar\joinrel\relbar}\limits^{*}_{#1}}}
\def\etra#1{\smash{\mathop{\mid\kern
-1pt\joinrel\relbar\joinrel\relbar}\limits_{#1}}}

\def\iff{\Leftrightarrow}
\def\Rw{\Rightarrow}
\def\oo{\overline}

\def\wt{\widetilde}
\def\wh{\widehat}

\def\im{\mbox{Im}}
\def\haus{\mbox{Haus}}

\def\N{\mathbb{N}}
\def\RR{\mathbb{R}}

\def\aut{\mbox{Aut}}
\def\inn{\mbox{Inn}}

\def\per{\mbox{Per}}

\def\ker{\mbox{Ker}}

\def\max{\mbox{max}}
\def\min{\mbox{min}}

\def\sup{\mbox{sup}}

\def\ZZ{\mathbb{Z}}

\def\p{\varphi}

\def\inv{^{-1}}

\def\bi{\begin{itemize}}
\def\ei{\end{itemize}}
\def\beq{\begin{equation}}
\def\eeq{\end{equation}}

\newtheorem{T}{Theorem}[section]
\newcommand{\bt}{\begin{T}}
\newcommand{\et}{\end{T}}
\newcommand{\ftd}{$\square$\end{T}}

\newtheorem{Proposition}[T]{Proposition}
\newcommand{\bp}{\begin{Proposition}}
\newcommand{\ep}{\end{Proposition}}
\newcommand{\fpd}{$\square$\end{Proposition}}

\newtheorem{Lemma}[T]{Lemma}
\newcommand{\bl}{\begin{Lemma}}
\newcommand{\el}{\end{Lemma}}
\newcommand{\fld}{$\square$\end{Lemma}}

\newtheorem{Corol}[T]{Corollary}
\newcommand{\bc}{\begin{Corol}}
\newcommand{\ec}{\end{Corol}}
\newcommand{\fcd}{$\square$\end{Corol}}

\newtheorem{Remark}[T]{Remark}
\newcommand{\br}{\begin{Remark}}
\newcommand{\er}{\end{Remark}}
\newcommand{\frd}{$\square$\end{Remark}}

\newtheorem{Example}[T]{Example}
\newcommand{\be}{\begin{Example}}
\newcommand{\ee}{\end{Example}}

\newtheorem{Problem}[T]{Problem}
\newcommand{\bq}{\begin{Problem}}
\newcommand{\eq}{\end{Problem}}

\newcommand{\proof}
   {\par\medbreak\noindent{\bf Proof}.\enspace}

\newcommand{\qed}{
$\Box$
\par\bigbreak}

\normalbaselines
\def\abstract#1{\par\bigskip
\begingroup\small
\baselineskip=12truept
\begin{center}ABSTRACT\end{center}
\par\medskip\par\noindent
\null\hfill\hbox{\vbox{\hsize=5truein\noindent#1}}
\hfill\null\par\endgroup\par}

\title{H\"older conditions for endomorphisms of hyperbolic groups}
\author{{\bf V\'\i tor Ara\'ujo and Pedro V. Silva}}

\date{\today}

\begin{document}
\maketitle

\begin{center}\small
2010 Mathematics Subject Classification: 20F67, 20E36, 51M10

\bigskip

Keywords: hyperbolic groups, virtually free groups, endomorphisms,
uniform continuity, H\"older conditions, Lipschitz conditions
\end{center}

\abstract{It is proved that an endomorphism $\p$ of a hyperbolic group
  $G$ satisfies a H\"older condition with respect to a
  visual metric if and only if $\p$ is virtually injective and $G\p$
  is a quasiconvex subgroup of $G$. If $G$ is virtually free or
  torsion-free co-hopfian, then 
$\p$ is uniformly continuous if and only if it satisfies a H\"older
condition if and only 
if it is virtually injective. Lipschitz conditions are discussed for
free group automorphisms.}

\section{Introduction}

The concept of boundary of a free group has been for a number of years
a major subject of research from geometric, topological, dynamical,
algebraic or combinatorial viewpoints. The boundary of $F_A$, denoted
by $\partial F_A$, can be 
defined as the set of all infinite reduced words on $\wt{A} = A \cup
A\inv$, but the topological (metric) structure is of utmost
importance. It can be defined through the {\em prefix metric}. Given $u,v
\in F_A$, let $u \wedge v$ denote the longest common prefix of $u$ and
$v$. An ultrametric $p_A:F_A \to F_A \to \RR_0^+$ is defined by
$$p_A(u,v) = \left\{
\begin{array}{ll}
2^{-|u\wedge v|}&\mbox{ if } u\neq v\\
0&\mbox{ otherwise}
\end{array}
\right.$$
The completion $(\wh{F_A},\wh{p_A})$ can be described as 
$$\wh{F_A} = F_A \cup \partial F_A,$$
and the metric $\wh{p_A}$ is nothing but the prefix metric defined
for finite and infinite reduced words altogether.

The theory of hyperbolic groups generalizes many aspects of free
groups, and we can endow the boundary of a hyperbolic group with a
metric structure proceeding analogously. This can be achieved with the
help of the Gromov product and the visual metrics
$\sigma_{p,\gamma}^A$. If $G = F_A$, $p = 1$ and $\gamma = \ln 2$,
then $\sigma_{p,\gamma}^A$ is precisely the prefix metric defined
above. 

The completion $(\wh{G},\wh{\sigma_{p,\gamma}^A})$ of
$(G,\sigma_{p,\gamma}^A)$ produces the boundary $\partial G = \wh{G}
\setminus G$ and its metric structure, which induces the Gromov
topology on $\partial G$. This same topology can be induced by any of
the visual metrics $d \in V^A(p,\gamma,T)$. These are the metrics
considered in this paper, and their extensions $\wh{d}$ to $\wh{G}$. 

Since the completion  $(\wh{G},\wh{d})$ is also
compact, the endomorphisms of $G$ which admit a continuous extension
to the boundary are precisely the uniformly continuous ones. It is thus
a natural problem to determine which endomorphisms of $G$ admit such a
continuous extension. It is well known that automorphisms do.

Uniform continuity is implied by a H\"older condition. A mapping
$\p:(X,d) \to (X',d')$ satisfies a {\em H\"older condition} of
exponent $r > 0$ if there exists a constant $K > 0$ such that
$$d'(x\p,y\p) \leq K(d(x,y))^r$$
for all $x,y \in X$. A H\"older condition of
exponent $1$ is a {\em Lipschitz condition}. Using an analogy with
notions from complexity theory in theoretical computer science (and
inverting the exponent), we may identify 
H\"older conditions with {\em polynomial complexity} and Lipschitz
conditions with {\em linear complexity}.

In this paper, we are interested mainly on H\"older conditions
for endomorphisms, with respect to visual metrics. Given the
exponential in the definition of the visual
metric, it is not surprising that this reduces to some Lispschitz type condition
involving Gromov products. As a preliminary result, we show that
all visual metrics on a hyperbolic group are H\"older equivalent.

In the main theorem of the paper (Theorem \ref{hpc}), we establish
several equivalent 
conditions for a nontrivial endomorphism of a
hyperbolic group to satisfy a H\"older condition. The most interesting
is undoubtedly the last 
one: $p$ must be virtually injective and $G\p$ must be a quasiconvex
subgroup of $G$. This second requirement may be removed if the group
is virtually free or torsion-free co-hopfian, when we also show that
all uniformly continuous 
endomorphisms satisfy a H\"older condition. The second author had proved
in \cite[Proposition 7.2]{Sil5} that a nontrivial endomorphism of a finitely
generated virtually free group is
uniformly continuous if and only if it is virtually injective.
We ignore whether this is also true for hyperbolic groups. 

We discuss also Lipschitz conditions for automorphisms. It is easy to see
that every inner automorphism of a hyperbolic group satisfies a
Lipschitz condition, but we have only succeeded on finding a precise
characterization in the free group case. With respect to the canonical basis,
Lipschitz conditions occur only for compositions of permutation
automorphisms with inner automorphisms. If we allow arbitrary finite
generating sets, we have only the inner automorphisms, and the same
happens if we allow arbitrary bases in rank $\geq 3$. In rank $2$, we
obtain an intermediate class of automorphisms. 

One of the motivations for our work is the possibility of defining 
new pseudometrics on $\aut(G)$ for every hyperbolic group $G$. Given a
virtually injective endomorphism $\p$ of $G$ and a visual metric $d$
on $G$, write 
$$||\p||_d = \ln(\inf\{ r \geq 1 \mid \p \mbox{ satisfies a H\"older
  condition of exponent }r\inv \mbox{ with respect to }d\}).$$
Since 
\beq
\label{sno}
||\p\psi||_d \leq ||\p||_d + ||\psi||_d
\eeq
for all virtually injective endomorphisms $\p,\psi$ of $G$, we call
$||\cdot||_d$ a seminorm. All inner automorphisms have seminorm 0.

Now we define a pseudometric $\oo{d}$ on $\aut(G)$ by 
$$\oo{d}(\p,\psi) = \max\{ ||\p\inv\psi||_d, ||\psi\inv\p||_d \}.$$
The inequality (\ref{sno}) implies the triangular inequality for $\oo{d}$.
This pseudometric is the object of ongoing work by the authors.

The paper is organized as follows. In Section 2 we present basic concepts
and notation for hyperbolic groups. Visual metrics and some of their
basic properties in connection with H\"older conditions are discussed in
Section 3. The main results of the paper, characterizing which
uniformly continuous endomorphisms satisfy a H\"older condition, are
presented in Section 4. Simplifications for the case of virtually free
or torsion-free co-hopfian hyperbolic groups are discussed in Section
5. In Section 6, we discuss Lipschitz conditions. Finally, some
open problems are proposed in Section 7.

\section{Hyperbolic groups}
\label{hg}

We present in this section well-known facts regarding hyperbolic
spaces and hyperbolic groups. The reader is referred to \cite{BH,GH} for
details.

A mapping $\p:(X,d) \to (X',d')$ between metric spaces is called
an {\em isometric embedding} if $d'(x\p,y\p) = d(x,y)$ for all $x,y
\in X$. A surjective isometric embedding is an
{\em isometry}.

A metric space $(X,d)$ is said to be {\em geodesic} if, for all $x,y
\in X$, there exists an isometric embedding 
$\xi:[0,s] \to X$ such that $0\xi = x$ and $s\xi = y$, where $[0,s]
\subset \RR$ is endowed with the usual metric of $\RR$. 
We call $\xi$ a {\em geodesic} of
$(X,d)$. We shall often call $\im(\xi)$ a geodesic as well. In this
second sense, we may use the notation $[x,y]$ to denote an arbitrary geodesic
connecting $x$ to $y$.
Note that a geodesic metric space is always (path) connected.

A {\em quasi-isometric embedding} of metric spaces is a mapping $\p:(X,d) \to
(X',d')$ such that there exist constants $\lambda \geq 1$ and $K \geq 0$
satisfying
$$\frac{1}{\lambda}d(x,y) -K \leq d'(x\p,y\p) \leq \lambda d(x,y) +K$$
for all $x,y \in X$. We may call it a
$(\lambda,K)$-quasi-isometric embedding if we want to stress the
constants. 
If in addition
$$\forall x' \in X'\; \exists x \in X:\; d'(x',x\p) \leq K,$$
we say that $\p$ is a {\em quasi-isometry}.

Two metric spaces $(X,d)$ and $(X',d')$ are said to be {\em
  quasi-isometric} if there exists a quasi-isometry $\p:(X,d) \to
(X',d')$. Quasi-isometry turns out to be an equivalence relation on
the class of metric spaces.
A {\em quasi-geodesic} of $(X,d)$ is a quasi-isometric
embedding 
$\xi:[0,s] \to X$ such that $0\xi = x$ and $s\xi = y$, where $[0,s]
\subset \RR$ is endowed with the usual metric of $\RR$. 

Let $(X,d)$ be a geodesic metric space.
Given $x_0,x_1,x_2 \in X$, a {\em geodesic triangle} $[[x_0,x_1,x_2]]$
is a collection of
three geodesics $[x_0,x_1]$, $[x_1,x_2]$ and $[x_2,x_0]$ in $X$. 

Given $\delta \geq 0$, we say
that $(X,d)$ is $\delta$-{\em hyperbolic} if
\beq
\label{rips}
\forall y \in [x_0,x_2] \hspace{.3cm} d(y,[x_0,x_1] \cup
[x_1,x_2]) \leq \delta
\eeq 
holds for every geodesic triangle $[[x_0,x_1,x_2]]$ in 
$X$. If this happens for some $\delta \geq 0$, we say
that $(X,d)$ is {\em 
  hyperbolic}. 

Given $Y,Z \subseteq X$ nonempty, the {\em Hausdorff distance} between
$Y$ and $Z$ is defined by
$$\haus(Y,Z) = \max\{ \sup_{y\in Y} d(y,Z),\; \sup_{z\in Z} d(z,Y) \}.$$

If $(X,d)$ is $\delta$-hyperbolic and $\lambda \geq 1$, $K \geq 0$, it
follows from \cite[Theorem 5.4.21]{GH}
that there exists a constant $R(\delta,\lambda,K)$,
depending only on $\delta,\lambda,K$, such that any geodesic
and $(\lambda,K)$-quasi-geodesic in $X$ having the same initial and
terminal points 
lie at Hausdorff distance $\leq R(\delta,\lambda,K)$ from each
other. This constant will be used in the proof of several results.

\medskip

Given a subset $A$ of a group $G$, we denote by $\langle A\rangle$ the
subgroup of $G$ generated by $A$. We assume throughout the paper that
generating sets are finite. 

Given $G = \langle A\rangle$, we write $\wt{A} =
A \cup A\inv$. The {\em Cayley
  graph} $\Gamma_A(G)$ has vertex set $G$ and edges of the form $g
\mapright{a} ga$ for all $g \in G$ and $a \in \wt{A}$.
The {\em geodesic metric} $d_A$ on $G$ is defined by taking $d_A(g,h)$
to be the length of the shortest path connecting $g$ to $h$ in $\Gamma_A(G)$. 

Since $\im(d_A) \subseteq \N$, then $(G,d_A)$ is not a geodesic metric
space. However, we can remedy that by embedding $(G,d_A)$
isometrically into the {\em geometric 
realization} $\oo{\Gamma}_A(G)$ of $\Gamma_A(G)$, when vertices become
points and edges 
become segments of length 1 in some (euclidean) space, intersections
being determined by adjacency only. With the obvious metric,
$\oo{\Gamma}_A(G)$ is a geodesic metric space, and the geometric 
realization is unique up to isometry. We denote also by $d_A$ the
induced metric on $\oo{\Gamma}_A(G)$.

We say that the group $G$ is {\em hyperbolic} if the geodesic metric
space $(\oo{\Gamma}_A(G),d_A)$ is
hyperbolic.

If $A'$ is an alternative finite generating set of $G$ and
\beq
\label{buildqg1}
N_{A,A'} = {\rm max}(\{ d_{A'}(1,a) \mid a \in A \} \cup \{ d_{A}(1,a') \mid a'
\in A' \}),
\eeq
it is immediate that
\beq
\label{gmpol33}
\frac{1}{N_{A,A'}}d_{A'}(g,h) \leq d_{A}(g,h) \leq N_{A,A'}d_{A'}(g,h)
\eeq
holds for all $g,h \in G$, hence
the identity mapping $(G,d_A) \to (G,d_{A'})$ is a
quasi-isometry. It follows easily that the concept of hyperbolic group
is independent from 
the finite generating set considered, but the hyperbolicity
constant $\delta$ may vary with the generating set. 

Condition (\ref{rips}), which became the most popular way of
defining hyperbolic group, is known as {\em Rips condition}. An
alternative approach is given by the concept of Gromov product, which
we now define. It can be defined for every metric space.

Given $g,h,p \in G$, we 
define 
$$(g|h)_p^A = \frac{1}{2}(d_A(p,g) + d_A(p,h)-d_A(g,h)).$$
We say that $(g|h)_p^A$ is the {\em Gromov product} ot $g$ and $h$,
taking $p$ as basepoint. 

The following result is well known:

\bp
\label{hyp}
The following
conditions are equivalent for a group $G = \langle A\rangle$:
\bi
\item[(i)] $G$ is hyperbolic;
\item[(ii)] there exists some $\delta \geq 0$ such that
\beq
\label{hyp1}
(g_0|g_2)^A_p \geq {\rm min}\{ (g_0|g_1)^A_p, (g_1|g_2)^A_p \} -\delta
\eeq
holds for all $g_0,g_1,g_2,p \in G$.
\ei
\ep

Let $H$ be a subgroup of a hyperbolic group $G = \langle A\rangle$ and
let $q \geq 0$. We
say that $H$ is $q$-{\em quasiconvex} with respect to $A$ if
$$\forall x \in [h,h']\hspace{.3cm} d_A(x,H) \leq q$$
holds for every geodesic $[h,h']$ in $\oo{\Gamma}_A(G)$ with endpoints in
$H$. We say that $H$ is {\em quasiconvex} if it is $q$-quasiconvex
for some $q \geq 0$.
Like most other properties in the theory of hyperbolic groups,
quasiconvex does not depend on the finite generating set considered
\cite[Section III.$\Gamma$.3]{BH}. 

A (finitely generated) subgroup of a hyperbolic group needs not be
hyperbolic, but a quasiconvex subgroup of a hyperbolic group is
always hyperbolic. The converse is not true in general.  
Quasiconvex subgroups occur quite frequently in the theory of
hyperbolic groups. In fact, non quasi convex subgroups are a 
relatively rare phenomenon, see \cite{Kap}.

We present next a model for the boundary of $G$.

Given a mapping $\p:\N \times \N \to \RR$, we write 
$$\displaystyle\varliminf_{i,j \to +\infty} (i,j)\p = \lim_{n \to
 +\infty}(\inf\{ (i,j)\p \mid i,j \geq n\}).$$

Fix a generating set $A$ for $G$ and $p \in G$.
We say that a sequence $(g_n)_n$ on $G$ is a {\em Gromov sequence} if
$$\displaystyle\varliminf_{i,j \to +\infty} (g_i|g_j)^A_p = +\infty.$$
This property is independent from both $A$ and $g$. Two Gromov
sequences $(g_n)_n$ and $(h_n)_n$ on $G$ are {\em equivalent} if 
$$\displaystyle\lim_{n \to +\infty} (g_n|h_n)^A_p = +\infty.$$
We denote by $[(g_n)_n]$ the equivalence class of the Gromov
sequence $(g_n)_n$. The set of all such equivalence classes is one of
the standard models for the boundary $\partial G$, and is adopted
in this paper.

We can identify $G$ with the set of all constant sequences $(g)_n$
on $G$, and consider
$$\wh{G} = \partial G \cup \{ \{ (g)_n \} \mid g \in G\}.$$
The Gromov product is extended to $\wh{G}$ by setting
$$(\alpha|\beta)^A_p = \displaystyle\sup\{ \varliminf_{i,j \to +\infty}
(g_i|h_j)^A_p \mid (g_n)_n \in \alpha,\; (h_n)_n \in \beta \}$$
for all $\alpha,\beta \in \wh{G}$.

\section{The visual metrics}

Let $G = \langle A\rangle$ be a hyperbolic group. Assuming that
$\Gamma_A(G)$ is $\delta$-hyperbolic, let $\gamma 
> 0$ be such that $\gamma\delta \leq \ln 2$.
Following Holopainen, Lang and
V\"ah\"akangas \cite{HLV}, we define 
$$\rho^A_{p,\gamma}(g,h) =
\left\{
\begin{array}{ll}
e^{-\gamma(g|h)_p^A}&\mbox{ if } g\neq h\\
0&\mbox{ otherwise}
\end{array}
\right.$$
for all $p,g,h \in G$. In general, $\rho^A_{p,\gamma}$ fails to be a
metric because of the
triangular inequality. Let
$$\sigma^A_{p,\gamma}(g,h) = 
\inf\{ \rho^A_{p,\gamma}(x_0,x_1) + \ldots + \rho^A_{p,\gamma}(x_{n-1},x_n) \mid 
n \geq 0,\; x_0 = g, \; x_n = h;\; x_1,\ldots,x_{n-1} \in G \}.$$
By \cite{HLV} (cf. also \cite{Fri,Vai}),
$\sigma^A_{p,\gamma}$ is a metric on
$G$ and the inequalities 
\beq
\label{ineq0}
\frac{1}{4}\rho^A_{p,\gamma}(g,h) \leq \sigma^A_{p,\gamma}(g,h) \leq
\rho^A_{p,\gamma}(g,h) 
\eeq
hold for all $g,h \in G$. 

The metric $\sigma^A_{p,\gamma}$ is an important example of a {\em
  visual metric}. Given $p \in G$, $\gamma > 0$ and $T \geq 1$, we
denote by $V^A(p,\gamma,T)$ the set of all metrics $d$ on $G$ such
that
\beq
\label{ineq}
\frac{1}{T}\rho^A_{p,\gamma}(g,h) \leq d(g,h) \leq
T\rho^A_{p,\gamma}(g,h)
\eeq
holds for all distinct $g,h \in G$. By (\ref{ineq0}), we have 
$$\sigma^A_{p,\gamma} \in V^A(p,\gamma,4)$$
whenever $\gamma\delta \leq \ln 2$. We shall refer to the metrics in some
$V^A(p,\gamma,T)$ as the {\em visual metrics} on $G$. 

Let $d \in V^A(p,\gamma,T)$ be a visual metric. In general, the metric
space $(G,d)$ is not 
complete. But its completion is essentially unique and also compact, and
can be obtained by 
adding to $G$ the elements of the boundary $\partial G$
\cite{BH,GH,HLV,Vai}. We denote it by $(\wh{G},\wh{d})$. It is well
known that $\wh{d}$ induces the {\em Gromov topology} on $\partial G$ 
\cite[Section III.H.3]{BH}. 

To understand the metric $\wh{d}$, we must consider
the extension of $\rho^A_{p,\gamma}$ to the boundary.
We define 
$$\hat{\rho}^A_{p,\gamma}(\alpha,\beta) =
\left\{
\begin{array}{ll}
e^{-\gamma(\alpha|\beta)_p^A}&\mbox{ if } \alpha\neq \beta\\
0&\mbox{ otherwise}
\end{array}
\right.$$
for all $\alpha,\beta \in \wh{G}$. By continuity, the inequalities 
\beq
\label{ineq2}
\frac{1}{T}\hat{\rho}^A_{p,\gamma}(\alpha,\beta) \leq
\wh{d}(\alpha,\beta) \leq
T\hat{\rho}^A_{p,\gamma}(\alpha,\beta) 
\eeq
hold for all $\alpha,\beta \in \wh{G}$ \cite[Section III.H.3]{BH}. 

It is widely known that uniform continuity of a mapping $\p:G \to G'$ of
hyperbolic groups determines the existence of a continuous
extension $\Phi:\wh{G} \to \wh{G'}$:

\bl
\label{ece}
Let $\p:G \to G'$ be a mapping of
hyperbolic groups and let $d$ and $d'$ be visual
metrics on $G$ and $G'$ respectively. Then the 
following conditions are equivalent:
\bi
\item[(i)] $\p$ is uniformly continuous with respect to $d$ and $d'$;
\item[(ii)] $\p$ admits a continuous extension $\Phi:(\wh{G},\wh{d})
  \to (\wh{G'},\wh{d'})$.
\ei
\el

Indeed, by a general topology result \cite[Section
XIV.6]{Dug}, every uniformly 
continuous mapping $\p:G \to G'$ admits a continuous extension to the
completions.

On the other hand, the completion $(\wh{G},\wh{d})$ is compact. Since every
continuous mapping with compact domain is uniformly continuous, it
follows that $\p$, being a restriction of a uniformly continuous
extension, is itself uniformly continuous.

We note also that the continuous extension is uniquely defined through
$$[(g_n)_n]\Phi = [(g_n\p)_n],$$
for every Gromov sequence $(g_n)_n$ on $G$.

A group is {\em virtually free} if it has a free subgroup of finite
index. Finitely generated virtually free groups constitute an
important subclass of hyperbolic groups. We should mention that the
second author developed 
in \cite{Sil5} a model for the boundary of such a group which allows
a huge simplification with respect to the general case.

\bl
\label{oldregp}
Let $G$ be a hyperbolic group and let $d \in V^A(p,\gamma,T)$, $d' 
\in V^{A'}(p',\gamma',T')$ be visual metrics on
$G$. Let $\p: (G,d) \to (G,d')$ be a mapping and let $P >
0$ and $Q \in \RR$ be constants such that
\beq
\label{oldregp1}
P(g\p|h\p)_{p'}^{A'} + Q \geq (g|h)_p^A
\eeq
holds for all $g,h \in G$.
Then $\p$ satisfies a H\"older condition of
exponent $\frac{\gamma'}{\gamma P}$.
\el

\proof
Let $g,h \in G$. We may assume that $g\p \neq h\p$, hence
$$\begin{array}{lll}
d'(g\p,h\p)&\leq&T'\rho^{A'}_{p',\gamma}(g\p,h\p) =
T'e^{-\gamma'(g\p|h\p)_{p'}^{A'}} \leq
T'e^{-\frac{\gamma'}{P}((g|h)_{p}^{A} -Q)}\\
&=&T'e^{\frac{\gamma' Q}{P}}e^{-\frac{\gamma'}{P}(g|h)_{p}^{A}} =
T'e^{\frac{\gamma'
    Q}{P}}(e^{-\gamma(g|h)_{p}^{A}})^{\frac{\gamma'}{\gamma P}} =
T'e^{\frac{\gamma'
    Q}{P}}(\rho^A_{p,\gamma}(g,h))^{\frac{\gamma'}{\gamma P}}\\
&\leq&T'e^{\frac{\gamma'
    Q}{P}}(Td(g,h))^{\frac{\gamma'}{\gamma P}}
= T'e^{\frac{\gamma'
    Q}{P}}T^{\frac{\gamma'}{\gamma P}}(d(g,h))^{\frac{\gamma'}{\gamma P}}
\end{array}$$
and we are done.
\qed

Replacing (\ref{ineq}) by (\ref{ineq2}), we may use the same argument
to prove the following:

\bl
\label{bold}
Let $G$ be a hyperbolic group and let $d \in V^A(p,\gamma,T)$, $d' 
\in V^{A'}(p',\gamma',T')$ be visual metrics on
$G$. Let $\Phi: (\wh{G},\wh{d}) \to (\wh{G},\wh{d'})$ be a mapping and
let $P > 0$ and $Q \in 
\RR$ be constants such that
$$P(\alpha\Phi|\beta\Phi)_{p'}^{A'} + Q \geq (\alpha|\beta)_p^A$$
holds for all $\alpha,\beta \in \wh{G}$.
Then $\Phi$ satisfies a H\"older condition of
exponent $\frac{\gamma'}{\gamma P}$. 
\el


An endomorphism $\p$ of a group $G$ is {\em trivial} if $G\p = 1$.
We show next that in the case of nontrivial endomorphisms, H\"older conditions
with respect to visual metrics 
are equivalent to Lispschitz conditions involving the Gromov product.

\bp
\label{regp}
Let $G$ be a hyperbolic group and let $d = \sigma^A_{p,\gamma}$, $d' =
\sigma^{A'}_{p',\gamma'}$ be visual metrics on
$G$. Let $\p: (G,d) \to (G,d')$ be a nontrivial homomorphism and let $r >
0$. Then the following conditions are equivalent:
\bi
\item[(i)] $\p$ satisfies a H\"older condition of
exponent $r$;
\item[(ii)] there exists a constant $Q \in \RR$ such that
\beq
\label{regp1}
\frac{\gamma'}{r\gamma}(g\p|h\p)_{p'}^{A'} + Q \geq (g|h)_p^A
\eeq
holds for all $g,h \in G$.
\ei
\ep

\proof
(i) $\Rw$ (ii). There 
exists a constant $K > 0$  such that
$$d'(g\p,h\p) \leq K(d(g,h))^r$$
for all $g,h \in G$. 

Assume first that $g\p \neq h\p$. Then $g \neq h$ and
$$\begin{array}{lll}
e^{-\gamma'(g\p|h\p)^{A'}_{p'}}&=&\rho^{A'}_{p',\gamma'}(g\p,h\p) \leq
T'd'(g\p,h\p) \leq T'K(d(g,h))^r \leq T'K(T\rho^{A}_{p,\gamma}(g,h))^r\\
&\leq&T'KT^r(\rho^{A}_{p,\gamma}(g,h))^r = T'KT^re^{-r\gamma(g|h)^{A}_{p}},
\end{array}$$
hence
$$-\gamma'(g\p|h\p)^{A'}_{p'} \leq \ln(T'KT^r)-r\gamma(g|h)^{A}_{p}$$
and so
\beq
\label{dif1}
\frac{\gamma'}{r\gamma}(g\p|h\p)_{p'}^{A'} +
\frac{\ln(T'KT^r)}{r\gamma} \geq (g|h)_p^A
\eeq
holds whenever $g\p \neq h\p$.

Now, since $\p$ is nontrivial, there exists some $a \in A$ such that $a\p
\neq 1$. We may assume that $d_A(1,a\p)$ is minimal.
We show that (\ref{regp1}) holds for 
$$Q = 1 + \frac{\ln(T'KT^r)}{r\gamma} + \frac{\gamma'}{r\gamma}d_{A'}(1,a\p).$$

In view of (\ref{dif1}), we may assume that $g\p = h\p$. On the one
hand, using (\ref{dif1}), we have
$$\begin{array}{lll}
(g\p|h\p)_{p'}^{A'}&=&\frac{1}{2}(d_{A'}(p',g\p)+d_{A'}(p',h\p)
-d_{A'}(g\p,h\p))\\ 
&\geq&\frac{1}{2}(d_{A'}(p',g\p)+d_{A'}(p',(ha)\p)
-d_{A'}(g\p,(ha)\p)-2d_{A'}(h\p,(ha)\p))\\
&=&(g\p|(ha)\p)_{p'}^{A'} -d_{A'}(1,a\p)\\
&\geq&\frac{r\gamma}{\gamma'}(g|ha)^A_p -\frac{\ln(T'KT^r)}{\gamma'}
-d_{A'}(1,a\p). 
\end{array}$$

On the other hand, $a\p \neq 1$ implies $a \neq 1$ and so
$$\begin{array}{lll}
(g|ha)_{p}^{A}&=&\frac{1}{2}(d_{A}(p,g)+d_{A}(p,ha)
-d_{A}(g,ha))\\ 
&\geq&\frac{1}{2}(d_{A}(p,g)+d_{A}(p,h)
-d_{A}(g,h)-2d_{A}(h,ha))\\
&=&(g|h)_{p}^{A} -1, 
\end{array}$$
hence
$$\begin{array}{lll}
(g\p|h\p)_{p'}^{A'}&\geq&
\frac{r\gamma}{\gamma'}(g|ha)^A_p -\frac{\ln(T'KT^r)}{\gamma'}
-d_{A'}(1,a\p) \geq \frac{r\gamma}{\gamma'}((g|h)^A_p-1)
-\frac{\ln(T'KT^r)}{\gamma'} -d_{A'}(1,a\p)\\
&=&\frac{r\gamma}{\gamma'}((g|h)_{p}^{A} -Q)
\end{array}$$
and so (\ref{regp1}) holds as required.

(ii) $\Rw$ (i). By Lemma \ref{oldregp}.
\qed

The next technical lemma illustrates an easy way of producing
quasi-geodesics. If $(X,d)$ is a geodesic metric space and
\beq
\label{gbp}
x_0 \mapright{} x_1 \mapright{} \ldots \mapright{} x_n
\eeq
is a path in $X$ such that each $x_{i-1} \mapright{} x_i$ is a
geodesic, then (\ref{gbp}) induces a canonical mapping $\xi:[0,s] \to
X$ such that $s = d(x_0,x_1) + \ldots + d(x_{n-1},x_n)$, $0\xi = x_0$
and $s\xi = x_n$.

\bl
\label{newqg}
Let $(X,d)$ be a geodesic metric space and let $\xi:[0,s] \to
X$ be the canonical mapping induced by
$$x_0 \mapright{} x_1 \mapright{} \ldots \mapright{} x_n,$$
where each $x_{i-1} \mapright{} x_i$ is a geodesic.
Let $P,L > 0$ and $Q \geq 0$ be such that 
\beq
\label{newqg2}
1 \leq d(x_{i-1},x_i) \leq L
\eeq
and
\beq
\label{newqg3}
Pd(x_i,x_j)+Q \geq |i-j|
\eeq
for all $i,j \in \{ 1,\ldots,n\}$. Then $\xi$ is a
$(\lambda,K)$-quasigeodesic for 
$$\lambda = {\rm max}\{ 1, LP\}, \quad K = {\rm max}\{2L,\frac{Q+1}{P}
+L\}.$$
\el

\proof
For $k = 0,\ldots,n$, let
$$s_k = \sum_{i=1}^k d(x_{i-1},x_i).$$
Clearly, $s_k\xi = x_k$.
It suffices to show that
\beq
\label{newqg1}
\frac{|u-v|}{LP}-\frac{Q+1}{P} -L \leq d(u\xi,v\xi) \leq |u-v| + 2L
\eeq
for all $u,v \in [0,s]$.

Indeed, it follows from (\ref{newqg2}) that there exist some $i,j \in
\{0, \ldots,n\}$ such that
$$d(u\xi,s_i\xi) = |u-s_i| \leq \frac{L}{2}, \quad
d(v\xi,s_j\xi) = |v-s_j| \leq \frac{L}{2}.$$
By symmetry, we may assume that $i \leq j$. Hence
$$\begin{array}{lll}
d(u\xi,v\xi)&\leq&d(s_i\xi,s_j\xi) + L = d(x_i,x_j) + L\\
&\leq&\sum_{\ell=i+1}^{j} d(x_{\ell-1},x_{\ell}) + L = 
\sum_{\ell=i+1}^{j} (s_{\ell}-s_{\ell-1}) + L\\
&=&s_j-s_i+L \leq |u-v| +2L.
\end{array}$$

On the other hand, (\ref{newqg2}) and (\ref{newqg3}) yield
$$\begin{array}{lll}
d(u\xi,v\xi)&\geq&d(s_i\xi,s_j\xi) -L = d(x_i,x_j) -L\\
&\geq&\frac{1}{P}|i-j| - \frac{Q}{P} -L \geq
\frac{1}{P}\sum_{\ell=i+1}^{j} \frac{s_{\ell}-s_{\ell-1}}{L} -
\frac{Q}{P} -L\\
&=&\frac{s_j-s_i}{LP}-\frac{Q}{P} -L \geq
\frac{|u-v|}{LP}-\frac{Q+1}{P} -L
\end{array}$$
and so (\ref{newqg1}) holds as required.
\qed

Two metrics $d$ and $d'$ on a set $X$ are {\em H\"older equivalent} if
the identity mappings 
$(X,d) \to (X,d')$ and $(X,d')  
\to (X,d)$ satisfy both a H\"older condition.

The following proposition is the finite version of the well-known analogue
result on the equivalence of the visual metrics on the boundary (see
\cite[Theorem 2.18]{KB}).

\bp
\label{gmpol}
All visual metrics on a given hyperbolic group are H\"older
  equivalent.
\ep

\proof
Let $G$ be a hyperbolic group and let $d,d'$ be visual metrics on
$G$.
Let $A,A'$ be finite generating sets of $G$ and assume that
$\oo{\Gamma}_A(G)$ (respectively $\oo{\Gamma}_{A'}(G)$) is $\delta$-hyperbolic
(respectively $\delta'$-hyperbolic). Let $p,p' \in G$.
In view of Proposition \ref{regp}, it suffices to show that there exist
constants $P > 0$ and $Q \geq 0$ such that
\beq
\label{newt1}
P(g|h)_{p'}^{A'} + Q \geq (g|h)_p^A
\eeq
holds for all $g,h \in G$. 

Let $N = N_{A,A'}$ be as in (\ref{buildqg1}) and let $R =
R(\delta,N^2,2N)$ be the constant introduced in Section \ref{hg}. Let
$[g,h]_{A}$ and $[g,h]_{A'}$ be geodesics in 
$\oo{\Gamma}_{A}(G)$ and $\oo{\Gamma}_{A'}(G)$, respectively. We claim that 
\beq
\label{gmpol2}
d_{A}(p,[g,h]_{A}) \leq Nd_{A'}(p',[g,h]_{A'}) +  N + d_{A}(p,p')+
R. 
\eeq

Assume that $[g,h]_{A'}$ is the path
$$g = g_0 \mapright{a'_1} g_1 \mapright{a'_2} \ldots \mapright{a'_n} g_n
= h$$
with $a'_1,\ldots,a'_n \in \wt{A'}$. Consider geodesics $g_{i-1}
\mapright{u_i} g_i$ in $\oo{\Gamma}_A(G)$ and let $\xi:[0,s] \to
(\oo{\Gamma}_{A}(G),d_{A})$ be the canonical mapping induced by the path
$$g = g_0 \mapright{u_1} g_1 \mapright{u_2} \ldots \mapright{u_n} g_n
= h.$$
Then $1 \leq d_A(g_{i-1},g_i) \leq N$ and in view of (\ref{gmpol33})
$$Nd_A(g_i,g_j) \geq d_{A'}(g_i,g_j) = |i-j|$$
holds for all $i,j \in \{ 1,\ldots,n\}$.
By Lemma \ref{newqg}, $\xi$ is a
$(N^2,2N)$-quasi-geodesic. Note that
$0\xi = g$, $s\xi = h$ and $[g,h]_{A'} \cap G \subseteq
\im(\xi)$. 

Now $$\haus([g,h]_A,\im(\xi)) \leq R(\delta,N^2,2N) = R,$$
hence
\beq
\label{gmpol4}
d_{A}(p,[g,h]_{A}) \leq d_{A}(p,\im(\xi)) + R \leq
d_{A}(p,p')+d_{A}(p',\im(\xi)) + R.
\eeq
On the other hand, we have
$$d_{A'}(p',[g,h]_{A'}) \geq d_{A'}(p',[g,h]_{A'} \cap G) -1$$
and 
$$d_{A}(p',\im(\xi)) \leq d_{A}(p',\im(\xi)\cap G) \leq
d_{A}(p',[g,h]_{A'} \cap G)$$
follows from $[g,h]_{A'} \cap G \subseteq
\im(\xi)$.  
In view of (\ref{gmpol33}), we get
$$d_{A'}(p',[g,h]_{A'}) \geq d_{A'}(p',[g,h]_{A'} \cap G) -1 \geq
\frac{1}{N}d_{A}(p',[g,h]_{A'} \cap G) -1 \geq
\frac{1}{N}d_{A}(p',\im(\xi))-1.$$
Together with (\ref{gmpol4}), this yields (\ref{gmpol2}).

By
\cite[Lemmas 2.9, 2.31 and 2.32]{Vai}, we have
\beq
\label{fgvf4}
(g|h)_p^A \leq d_A(p,[g,h]_A) \leq (g|h)_p^A +2\delta.
\eeq

Together with (\ref{gmpol2}), this yields
$$(g|h)_p^A \leq Nd_{A'}(p',[g,h]_{A'}) +  N + d_{A}(p,p')+
R.$$
Applying (\ref{fgvf4}) to $d_{A'}(p',[g,h]_{A'})$, we obtain
$$(g|h)_p^A \leq N(g|h)_{p'}^{A'} + 2N\delta' +  N + d_{A}(p,p')+
R,$$
hence (\ref{newt1}) holds for $P = N$ and $Q =
(2\delta'+1)N + d_{A}(p,p')+ R$.
\qed

\section{Endomorphisms of hyperbolic groups}

An endomorphism $\p$ of $G$ is {\em virtually
  injective} if its kernel is finite. This is a necessary condition
for uniform continuity:

\bl
\label{ucv}
Let $G$ be a hyperbolic group endowed with a visual
metric $d$. Let $\p$ be a uniformly continuous nontrivial endomorphism
of $G$. Then $\p$ is virtually injective.
\el

\proof
Assume that $d \in
V^A(p,\gamma,T)$. Fix $g \in G \setminus \ker(\p)$. Let $\varepsilon =
d(1,g\p) > 0$ and let 
$\delta > 0$ be such that
$$\forall x,y \in G \; (d(x,y) < \delta \Rw d(x\p,y\p) <
\varepsilon).$$
For every $h \in \ker(\p)$, we have $d(h\p,(hg)\p) = d(1,g\p) =
\varepsilon$, hence $d(h,hg) \geq \delta$. By (\ref{ineq}), we get
$$e^{-\gamma(h|hg)_p} = \rho^A_{p,\gamma}(h,hg) \geq
\frac{1}{T}d(h,hg) \geq \frac{\delta}{T}$$
and so
$$(h|hg)_p \leq -\frac{\ln\frac{\delta}{T}}{\gamma}.$$
It follows that
$$\begin{array}{lll}
d_A(p,h)&\leq&\frac{1}{2}(d_A(p,h) + d_A(p,hg) -d_A(h,hg) + 2d_A(h,hg))
= (h|hg)_p + d_A(h,hg)\\
&=&(h|hg)_p + d_A(1,g) \leq
-\frac{\ln\frac{\delta}{T}}{\gamma} + d_A(1,g).
\end{array}$$
Since $A$ is finite, then $\Gamma_A(G)$ is locally finite, i.e. every
ball is finite. Therefore $\ker(\p)$ is finite and
$\p$ is virtually injective.
\qed

We need also the following result:

\bp
\label{eqq}
Let $\p$ be a nontrivial endomorphism of a hyperbolic group $G$ with
continuous extension $\Phi:\wh{G} \to \wh{G}$. Let $d$
be a visual metric on $G$. Then the  
following conditions are equivalent:
\bi
\item[(i)] $\p$ satisfies a H\"older condition of exponent $r$  with
  respect to $d$; 
\item[(ii)] $\Phi$ satisfies a H\"older condition of exponent $r$
  with respect to $\wh{d}$. 
\ei
\ep

\proof
(i) $\Rw$ (ii). Let $d \in V^A(p,\gamma,T)$.
By Proposition \ref{regp}, there exists a constant $Q
\in \RR$ such that
\beq
\label{eqq1}
\frac{1}{r}(g\p|h\p)_{p}^{A} + Q \geq (g|h)_p^A
\eeq
holds for all $g,h \in G$. We show that 
\beq
\label{eqq2}
\frac{1}{r}(\alpha\Phi|\beta\Phi)_{p}^{A} + Q \geq (\alpha|\beta)_p^A
\eeq
for all $\alpha,\beta \in \wh{G}$. 

Let $(g_n)_n \in \alpha$ and $(h_n)_n \in \beta$. For every $n \in
\N$, we have by (\ref{eqq1})
$$\inf\{ (g_i\p|h_j\p)_p^A \mid i,j \geq n\}
\geq r\cdot\inf\{ (g_i|h_j)_p^A \mid i,j \geq n\} -rQ.$$
Hence
$$\varliminf_{i,j \to +\infty} (g_i\p|h_j\p)_p^A
\geq r\varliminf_{i,j \to +\infty} (g_i|h_j)_p^A
-rQ.$$
It follows that
$$\sup\{ \varliminf_{i,j \to +\infty} (g_i\p|h_j\p)_p^A \mid 
(g_n)_n \in \alpha, (h_n)_n \in \beta \}
\geq r(\alpha|\beta)^A_p -rQ.$$
Since $(g_n)_n \in \alpha$ implies $(g_n\p)_n \in \alpha\Phi$ and 
$(h_n)_n \in \beta$ implies $(h_n\p)_n \in \beta\Phi$, it follows that
$$(\alpha\Phi|\beta\Phi)^A_p \geq r(\alpha|\beta)^A_p
-rQ.$$
By Lemma \ref{bold}, 
$\Phi$ satisfies a H\"older condition of exponent $r$  with respect to $\wh{d}$.
 
(ii) $\Rw$ (i). Immediate since $\Phi$ is an extension of $\p$ and
$\wh{d}$ is an extension of $d$.
\qed

In the main result of the paper, we characterize the uniformly
continuous endomorphisms which satisfy a H\"older condition:

\bt
\label{hpc}
Let $\p$ be a nontrivial endomorphism of a hyperbolic group $G$ and
let $d \in V^A(p,\gamma,T)$ be a visual metric on $G$. Then the 
following conditions are equivalent:
\bi
\item[(i)] $\p$ satisfies a H\"older condition with respect to $d$;
\item[(ii)] $\p$ admits an extension to $\wh{G}$ satisfying
  a H\"older condition with respect to $\wh{d}$; 
\item[(iii)] there exist constants $P > 0$ and $Q \in \RR$ such that
\beq
\label{hpc99}
P(g\p|h\p)_p^A + Q \geq (g|h)_p^A
\eeq
for all $g,h \in G$;
\item[(iv)] $\p$ is a quasi-isometric embedding of $(G,d_A)$ into itself;
\item[(v)] $\p$ is virtually injective and $G\p$ is a quasiconvex
  subgroup of $G$.
\ei
\et

\proof
(i) $\iff$ (ii). By Lemma \ref{ece} and Proposition \ref{eqq}.

(i) $\iff$ (iii). By Proposition \ref{regp}.

(i) $\Rw$ (v). By Lemma \ref{ucv}, $\p$ is virtually injective.
In view of Proposition \ref{gmpol}, we may assume that $p = 1$. Since
(i) implies (iii), there exist constants $P > 0$ and $Q \in \RR$ such that
$$P(g\p|g\p)_1^A + Q \geq (g|g)_1^A$$
for every $g \in G$, which is equivalent to
$$Pd_A(1,g\p) + Q \geq d_A(1,g).$$



Since $d_A(g,h) = d_A(1,g\inv h)$ and $d_A(g\p,h\p) = d_A(1,(g\inv
h)\p)$, we immediately get
\beq
\label{hpc5}
Pd_A(g\p,h\p) + Q \geq d_A(g,h)
\eeq
for all $g,h \in G$.

Let 
$$M_{\p} = \max\{ d_A(1,a\p) \mid a \in A\}.$$
We show now that $G\p$ is quasiconvex.

Let $g,h \in G$ and let $g' \mapright{w} h'$ have minimal length among
all the paths in $\oo{\Gamma}_A(G)$ such that $g'\p = g\p$ and $h'\p =
h\p$. In particular, $g' \mapright{w} h'$ is a geodesic. 
Assume that $w = a_1\ldots a_n$ with $a_i \in \wt{A}$. For $i = 0,
\ldots,n$, write $w_i = a_1\ldots a_i$ and let $(g'w_{i-1})\p
\mapright{} (g'w_{i})\p$ be a geodesic. Let $\xi:[0,s] \to
\oo{\Gamma}_A(G)$ be the canonical mapping induced by the path
$$g\p = g'\p = (g'w_{0})\p \mapright{} (g'w_{1})\p \mapright{} 
\ldots \mapright{} (g'w_{n})\p = h'\p = h\p.$$

Suppose that $a_i\p = 1$ for some $i$. Let $w' = w_{i-1}a_{i+1}\ldots
a_n$. Since $w'\p = w\p$, we have $(g'w')\p = (g'w)\p = h'\p = h\p$,
contradicting the minimality of $w$. Hence $a_i\p \neq 1$ for every
$i$. Since
$|a_i\p| \leq M_{\p}$, we get
$$1 \leq d_A((g'w_{i-1})\p,(g'w_{i})\p) \leq M_{\p}.$$
Assume that $0 \leq i \leq j \leq n$. By (\ref{hpc5}), we have
$$Pd_A((g'w_{i})\p,(g'w_{j})\p) + Q \geq 
d_A(g'w_{i},g'w_{j}) = d_A(1,a_{i+1}\ldots a_j).$$
Since $a_{i+1}\ldots a_j$ is a factor of a geodesic, it is itself a
geodesic and so $d_A(1,a_{i+1}\ldots a_j) = j-i$. Thus
$$Pd_A((g'w_{i})\p,(g'w_{j})\p) + Q \geq |j-i|$$
and it follows from Lemma \ref{newqg} that $\xi$ is a
$(\lambda,K)$-quasigeodesic for
$$\lambda = {\rm max}\{ 1, M_{\p}P\}, \quad K = {\rm max}\{2M_{\p},\frac{Q+1}{P}
+M_{\p}\}.$$

Let $R = R(\delta,\lambda,K)$. Let $[g\p,h\p]$ be a geodesic in
$\oo{\Gamma}_A(G)$ and let $x \in [g\p,h\p]$. Then 
$$\haus([g\p,h\p],\im(\xi)) \leq R$$
and every point in $\im(\xi)$ is at distance at most $M_{\p}$ from an
element of $G\p$, hence 
$$d_A(x, G\p) \leq R+M_{\p}$$
and so $G\p$ is $(R+M_{\p})$-quasiconvex.

(v) $\Rw$ (iv). Let $K = \ker(\p) \unlhd G$. Let $\pi:G \to G/K$ be
the canonical projection and let $\iota:G\p \to G$ be inclusion. Then
there exists an isomorphism $\oo{\p}: G/K \to G\p$ such that the
diagram
$$\xymatrix{
G \ar[r]^{\p} \ar[d]_{\pi} & G\\
G/K \ar[r]_{\oo{\p}} & G\p \ar[u]_{\iota}
}$$
commutes. Since the composition of quasi-isometric embeddings is still
a quasi-isometric embedding, it suffices to show that each one of the
homomorphisms $\pi,\oo{\p},\iota$ is a quasi-isometric embedding when
we consider a geodesic metric in each of the groups (it does not
matter which since the identity $(H,d_A) \to (H,d_B)$  is a
quasi-isometry whenever $H = \langle A \rangle = \langle B \rangle$ by
(\ref{gmpol33})). 

Let 
$$L = \max\{ d_A(1,x) \mid x \in K\}.$$
Let $g,h \in G$. We claim that
\beq
\label{fq}
d_A(g,h)-L \leq d_{A\pi}(g\pi,h\pi) \leq d_A(g,h).
\eeq

Since $h = ga_1\ldots a_n$ implies $h\pi = (ga_1\ldots a_n)\pi$ for
all $a_1,\ldots, a_n \in \wt{A}$, we have $d_{A\pi}(g\pi,h\pi) \leq
d_A(g,h)$.

Write $h\pi = (gw)\pi$, where $w$ is a word on $\wt{A}$ of minimum
length. Then $h = gwx$ for some $x \in K$ and so
$$d_A(g,h) \leq d_A(g,gw) + d_A(gw,gwx) = 
d_A(1,w) + d_A(1,x) \leq |w| + L.$$
By minimality of $w$, we have actually $|w| = d_{A\pi}(g\pi,h\pi)$ and
thus (\ref{fq}) holds.

Now $\oo{\p}$ is an isomorphism and $A\pi\oo{\p} = A\p$, hence
$$d_{A\p}(g\pi\oo{\p},h\pi\oo{\p}) = d_{A\pi}(g\pi,h\pi)$$
for all $g,h \in G$ and so $\oo{\p}:(G/K,d_{A\pi}) \to (G\p,d_{A\p})$
is actually an isometry.

Assume that $G\p$ is $q$-quasi convex with respect to $A$. Let 
$$B = \{ h \in G\p \mid d_A(1,h) \leq 2q+1\}.$$
Then $B$ is a finite generating set of $G\p$ and
$\iota:(G\p,d_{B}) \to (G,d_A)$ is a quasi-isometric embedding
\cite[Section III.$\Gamma$.3]{BH}. 
Therefore all three homomorphisms $\pi,\oo{\p},\iota$ are
quasi-isometric embeddings and so is their composition $\p$. 

(iv) $\Rw$ (iii).
Clearly, $\p$ can be extended to a quasi-isometric embedding $\oo{\p}$
of $(\oo{\Gamma}_A(G),d_A)$ into itself. Assume that $\oo{\Gamma}_A(G)$
is $\delta$-hyperbolic.
Let $\lambda \geq 1$ and $K \geq 0$ be constants such that
\beq
\label{hqe1}
\frac{1}{\lambda}d_A(x,y) -K \leq d_A(x\oo{\p},y\oo{\p}) \leq \lambda d_A(x,y)
+K
\eeq
holds for all $x,y \in \oo{\Gamma}_A(G)$. Write $R =
R(\delta,\lambda,K)$. We prove that
\beq
\label{hqe2}
\lambda(g\p|h\p)^A_p + \delta +\lambda(\lambda\delta + \frac{3K}{2} +
3R +d_A(p,p\p)) \geq (g|h)^A_p
\eeq
holds for all $g,h \in G$. 

Let $g,h \in G$. Consider a geodesic triangle $[[p,g,h]]$ with
geodesics $[p,g]$, $[g,h]$ and $[p,h]$. Let 
$$X = \{ x \in [g,h]: d_A(x,[p,g]) \leq \delta\}, \quad
Y = \{ y \in [g,h]: d_A(y,[p,h]) \leq \delta\}.$$
It is immediate that $X$ and $Y$ are both closed and nonempty.
Since $X \cup Y = [g,h]$ is obviously connected, it follows that $X
\cap Y \neq \emptyset$. Let $x \in X \cap Y$ and take $g' \in [p,g]$
and $h' \in [p,h]$ such that $d_A(x,g'), d_A(x,h') \leq \delta$.

If $\xi:[0,s] \to [p,g]$ is our geodesic, let $\xi' = \xi\oo{\p}$. For
all $i,j \in [0,s]$, (\ref{hqe1}) yields
$$d_A(i\xi'-j\xi') = d_A(i\xi\oo{\p}-j\xi\oo{\p}) \leq
\lambda|i\xi-j\xi|+K = \lambda|i-j|+K.$$
Similarly, 
$$d_A(i\xi'-j\xi') \geq \frac{1}{\lambda}|i-j|-K$$
and so $\xi'$ is a
$(\lambda,K)$-quasi-geodesic from $[0,s]$ to $\oo{\Gamma}_A(G)$ such
that $0\xi' = p\p$, $s\xi' = g\p$ and $g'\p \in \im(\xi')$. Fix a geodesic
$[p\p,g\p]$. Then  
$\haus([p\p,g\p],\im(\xi')]) \leq R$, hence there exists some $g'' \in
[p\p,g\p]$ such that $d_A(g'',g'\p) \leq R$. Similarly,
fix geodesics $[p\p,h\p]$ and
$[g\p,h\p]$. Then
there exist some $h'' \in
[p\p,h\p]$ and $x' \in [g\p,h\p]$ such that $d_A(h'',h'\p), d_A(x',x\p) \leq R$.

We claim that
\beq
\label{hqe3}
d_A(g'',g\p)-d_A(g\p,x') \geq -\lambda\delta - K - 2R.
\eeq
Indeed, in view of (\ref{hqe1}), we have
$$\begin{array}{lll}
d_A(g'',g\p)-d_A(g\p,x')&\geq&-d_A(g'',x') \geq -d_A(g'',g'\p) -
d_A(g'\p,x\p) - d_A(x\p,x')\\
&\geq&-\lambda d_A(g',x) -K -2R \geq -\lambda\delta -K-2R.
\end{array}$$
Similarly,
\beq
\label{hqe4}
d_A(h'',h\p)-d_A(h\p,x') \geq -\lambda\delta - K - 2R.
\eeq
Now (\ref{hqe1}), (\ref{hqe3}) and
(\ref{hqe4}) combined yield 
$$\begin{array}{lll}
(g\p|h\p)^A_{p\p}&=&\frac{1}{2}(d_A(p\p,g\p)+d_A(p\p,h\p)-d_A(g\p,h\p))\\
&=&\frac{1}{2}(d_A(p\p,g'')+
d_A(g'',g\p)+d_A(p\p,h'')+d_A(h'',h\p)-d_A(g\p,x')-d_A(x',h\p))\\ 
&=&\frac{1}{2}(d_A(p\p,g'')+d_A(p\p,h'')
+d_A(g'',g\p)-d_A(g\p,x')+d_A(h'',h\p)-d_A(x',h\p))\\ 
&\geq&\frac{1}{2}(d_A(p\p,g'\p)+d_A(p\p,h'\p)
+d_A(g'',g\p)-d_A(g\p,x')+d_A(h'',h\p)-d_A(x',h\p))-R\\
&\geq&\frac{1}{2}(d_A(p\p,g'\p)+d_A(p\p,h'\p)) - \lambda\delta - K - 3R\\
&\geq&\frac{1}{2\lambda}(d_A(p,g')+d_A(p,h')) - \lambda\delta -
\frac{3K}{2} - 3R. 
\end{array}$$
It follows that
$$\begin{array}{lll}
(g\p|h\p)^A_{p}&=&\frac{1}{2}(d_A(p,g\p)+d_A(p,h\p)-d_A(g\p,h\p))\\
&\geq&\frac{1}{2}(d_A(p\p,g\p)+d_A(p\p,h\p)-2d_A(p,p\p)-d_A(g\p,h\p))\\
&=&(g\p|h\p)^A_{p\p}-d_A(p,p\p)\\
&\geq&\frac{1}{2\lambda}(d_A(p,g')+d_A(p,h')) - \lambda\delta
- \frac{3K}{2} - 3R -d_A(p,p\p).
\end{array}$$

On the other hand,
$$\begin{array}{lll}
(g|h)^A_p&=&\frac{1}{2}(d_A(p,g)+d_A(p,h)-d_A(g,h))\\
&=&\frac{1}{2}(d_A(p,g')+
d_A(g',g)+d_A(p,h')+d_A(h',h)-d_A(g,x)-d_A(x,h))\\ 
&=&\frac{1}{2}(d_A(p,g')+d_A(p,h')
+d_A(g',g)-d_A(g,x) +d_A(h',h)-d_A(x,h))\\
&\leq&\frac{1}{2}(d_A(p,g')+d_A(p,h')
+d_A(g',x) +d_A(h',x))\\
&\leq&\frac{1}{2}(d_A(p,g')+d_A(p,h'))+\delta
\end{array}$$
and combining the previous ineqialities we get
$$\begin{array}{lll}
(g\p|h\p)^A_p&\geq&\frac{1}{2\lambda}(d_A(p,g'\p)+d_A(p,h'\p)) -
\lambda\delta - \frac{3K}{2} - 3R -d_A(p,p\p)\\
&\geq&\frac{1}{\lambda}(g|h)^A_p - \frac{\delta}{\lambda}
-\lambda\delta - \frac{3K}{2} - 3R -d_A(p,p\p).
\end{array}$$
Therefore (\ref{hqe2}) holds and we are done. 
\qed

\section{Simplifications}

Under which circumstances does every uniformly continuous endomorphism
satisfy a H\"older condition? We have the following remark:

\bl
\label{ucho}
Let $G$ be a hyperbolic
group and let $d \in
V^A(p,\gamma,T)$ be a visual metric on $G$. Then the  
following conditions are equivalent:
\bi
\item[(i)] every uniformly continuous endomorphism of $G$ satisfies a
  H\"older condition with respect to $d$; 
\item[(ii)] $G\p$ is a quasiconvex subgroup of $G$ for every
  endomorphism of $G$ uniformly continuous with respect to $d$.
\ei
\el

\proof
(i) $\Rw$ (ii). Both conditions hold for the trivial endomorphism.
For nontrivial endomorphisms we use Theorem \ref{hpc}. 

(ii) $\Rw$ (i). Let $\p$ be a nontrivial endomorphism of $G$, uniformly
continuous with respect to $d$. By Lemma \ref{ucv}, $\p$ is virtually
injective. Now we apply Theorem \ref{hpc}.
\qed

Now we get a simplified version of Theorem \ref{hpc} for
virtually free groups:

\bc
\label{vpc}
Let $\p$ be a nontrivial endomorphism of a finitely generated
virtually free group $G$ and let $d \in V^A(p,\gamma,T)$ be a visual
metric on $G$. Then the 
following conditions are equivalent:
\bi
\item[(i)] $\p$ is uniformly continuous with respect to $d$;
\item[(ii)] $\p$ satisfies a H\"older condition with respect to $d$; 
\item[(iii)] $\p$ admits a continuous extension to the completion
  $(\wh{G},\wh{d})$;
\item[(iv)] $\p$ admits an extension to $\wh{G}$ satisfying
  a H\"older condition with respect to $\wh{d}$; 
\item[(v)] there exist constants $P > 0$ and $Q \in \RR$ such that
$$P(g\p|h\p)_p^A + Q \geq (g|h)_p^A$$
for all $g,h \in G$;
\item[(vi)] $\p$ is a quasi-isometric embedding
of $(G,d_A)$ into itself;
\item[(vii)] $\p$ is virtually injective.
\ei
\ec

\proof
By \cite[Corollary 4.2]{AS}, every subgroup of a finitely generated
virtually free group 
is quasiconvex, hence the equivalences (ii) $\iff$ (iv) $\iff$ (v)
$\iff$ (vi) $\iff$ (vii) 
follow from Theorem \ref{hpc}. Now (ii) $\Rw$ (i) holds trivially, 
(i) $\Rw$ (vii) follows from Lemma \ref{ucv}, and (i) $\iff$ (iii)
follows from Lemma \ref{ece}. 
\qed

We note that the equivalence (i) $\iff$ (iii) $\iff$ (vii) had been
proved by the second author in \cite{Sil5}.

If $\p$ is not an endomorphism, then Corollary \ref{vpc} does not
hold, even for an infinite cyclic group. For $x \in \RR$, we denote by
$\lfloor x \rfloor$ the greatest integer $n \leq x$.

\be
\label{noten}
Let $d$ be the prefix metric on $F_{\{ a \}}$ and let $\p:F_{\{ a \}}
\to F_{\{ a \}}$ be defined by 
$$a^n\p = \left\{
\begin{array}{ll}
a^{\lfloor \sqrt{n} \rfloor}&\mbox{ if } n \geq 0\\
a^n&\mbox{ otherwise}
\end{array}
\right.$$
Then $\p$ is uniformly continuous but satisfies no H\"older condition 
with respect to $d$.
\ee

Indeed, it is a simple exercise to show that 
$$\forall \varepsilon > 0\; \forall m,n \in \ZZ\; 
(d(a^m,a^n) < \min\{ \frac{1}{2},
\frac{1}{2}\varepsilon^{2-\log_2\varepsilon} \}
 \Rw d(a^m\p,a^n\p) < \varepsilon),$$
hence $\p$ is uniformly continuous with respect to $d$.

However, for all $r,K > 0$, we have 
$$rm - \lfloor \sqrt{m} \rfloor > \log_2K$$
for $m$ large enough. It is easy to check that, for such $m$ and $n >
(\sqrt{m} +1)^2$, we have 
$$d(a^m\p,a^n\p) > K(d(a^m,a^n))^r,$$
thus $\p$ satisfies no H\"older condition with respect to $d$.

\bigskip

We would like to extend the previous discussion to the class of all
hyperbolic groups, but the panorama in not so clear. But we are able
to produce another result on the line of Corollary \ref{vpc}. Werecall
that a group $G$ is said to be {\em co-hopfian} if every monomorphism
of $G$ is an automorphism.

\bc
\label{tfch}
Let $\p$ be a nontrivial endomorphism of a torsion-free co-hopfian
hyperbolic group $G$ and let $d \in V^A(p,\gamma,T)$ be a visual
metric on $G$. Then the 
following conditions are equivalent:
\bi
\item[(i)] $\p$ is uniformly continuous with respect to $d$;
\item[(ii)] $\p$ satisfies a H\"older condition with respect to $d$; 
\item[(iii)] $\p$ admits a continuous extension to the completion
  $(\wh{G},\wh{d})$;
\item[(iv)] $\p$ admits an extension to $\wh{G}$ satisfying
  a H\"older condition with respect to $\wh{d}$; 
\item[(v)] there exist constants $P > 0$ and $Q \in \RR$ such that
$$P(g\p|h\p)_p^A + Q \geq (g|h)_p^A$$
for all $g,h \in G$;
\item[(vi)] $\p$ is a quasi-isometric embedding
of $(G,d_A)$ into itself;
\item[(vii)] $\p$ is virtually injective;
\item[(viii)] $\p$ is an automorphism.
\ei
\ec

\proof
(vii) $\iff$ (viii). Since $G$ is torsion-free, every virtually
injective endomorphism is a monomorphism. Then we use the fact that
$G$ is co-hopfian.

Now the equivalences (ii) $\iff$ (iv) $\iff$ (v)
$\iff$ (vi) $\iff$ (vii) 
follow from Theorem \ref{hpc}.

Finally, (ii) $\Rw$ (i) holds trivially, 
(i) $\Rw$ (vii) follows from Lemma \ref{ucv}, and (i) $\iff$ (iii)
follows from Lemma \ref{ece}. 
\qed

Examples of torsion-free co-hopfian hyperbolic groups have been
provided by Rips and Sela \cite{RSel} and Sela \cite{Sel},
namely:
\bi
\item
non-elementary torsion-free hyperbolic groups which admit no
nontrivial cyclic splittings \cite{RSel};
\item
non-elementary torsion-free freely indecomposable hyperbolic groups
\cite{Sel}.
\ei 

However, many torsion-free hyperbolic groups fail to be co-hopfian,
such as infinite cyclic groups. For more interesting examples, see
\cite{KW}.

In the case of torsion-free groups, it would be enough of course to
prevent the existence of monomorphisms with non quasiconvex image. The
following example shows that such a situation cannot always be avoided.

\be
\label{nqc}
There exists a torsion-free hyperbolic group $G$ having a non
quasiconvex subgroup isomorphic to $G$.
\ee

Indeed, let $a,b,t$ be distinct letters and write $A = \{ a,b\}$ and $B = \{
a,b,t\}$. We fix words $u = abab^2\ldots
ab^{20}$ and $v = baba^2\ldots
ba^{20}$. Let $H$ be the group defined by the presentation 
\beq
\label{piece}
\langle B \mid t\inv a tu, t\inv btv\rangle.
\eeq
Let $R$ denote the set of all cyclic conjugates of the two (cyclically
reduced) relators
and their inverses. A {\em piece} of (\ref{piece}) is a maximal common
prefix of two distinct elements of $R$. It is easy to see that the
longest pieces of (\ref{piece}) are $b^{18}ab^{19}$, $a^{18}ab^{19}$
and their inverses and have therefore length $38$. On the oher hand,
the length of each relator is $3+20+\frac{20(20+1)}{2} = 233$. Since
$38 < \frac{1}{6}233$, the presentation (\ref{piece}) satisfies the
small cancellation condition $C'(\frac{1}{6})$. Now it follows from a
theorem of Gromov \cite{Gro03} that $H$ is hyperbolic.

Let $F$ denote the subgroup of $H$ generated by $a,b$. The subgroup
$K$ of $F_A$ generated by $u$ and $v$ cannot have rank 1 since $uv
\neq vu$ in $F_A$. Since $F_A$ is hopfian, it follows that $K$ is free
on $\{ u,v\}$. Hence $H$ is an HNN extension of $F_A$ and so there is
a canonical isomorphism $F_A \to F$. Thus $F$ is a free subgroup of
$H$ with basis $A$. Moreover, since the finite order elements of the
HNN extension $H$ must be conjugates of the finite order elements of
$F$ (see \cite{LS}), then $H$ is torsion-free. 

Since $F$ is a normal subgroup of $H$, we have $t^{-n}at^n \in F$ for
every $n \geq 0$. Consider the geodesic metrics $d_A$ and $d_B$ on $F$
and $H$. We have
\beq
\label{distort1}
d_B(1,t^{-n}at^n) \leq 2n+1
\eeq
for every $n \geq 0$. We prove that
\beq
\label{distort2}
d_A(1,t^{-n}at^n) = 230^n
\eeq
by induction on $n$. The case $n = 0$ being trivial, assume that $n
\geq 1$ and $d_A(1,t^{-(n-1)}at^{n-1}) = 230^{n-1}$. It follows that
there exist $m = 230^{n-1}$ letters $c_1, \ldots, c_m \in \wt{A}$ such
that $t^{-(n-1)}at^{n-1} = c_1\ldots c_m$ in reduced form. Hence
$$t^{-n}at^n = (t\inv c_1t)\ldots (t\inv c_mt).$$
We have $t\inv c_it \in \{ u,v,u\inv,v\inv\}$ for $i =
1,\ldots,m$. Moreover, since $u = a\ldots b$ and $v = b\ldots a$, and
$c_1\ldots c_m$ is reduced, there is no reduction between the reduced
forms over $\wt{A}$ of two consecutive $t\inv c_it$. Hence the length
of the reduced form of $t^{-n}at^n$ over $\wt{A}$ is $230m =
230^n$. Since $F$ is free on $A$, we get (\ref{distort2}). 

Now it follows from (\ref{distort1}) and (\ref{distort2}) that the
embedding $(F,d_A) \to (H,d_B)$ is not a quasi-isometric embedding and
so $F$ is not an {\em undistorted} subgroup of $H$. By a theorem of
Short \cite{Sho} (see also \cite[Lemma $\Gamma$.3.5]{BH}), $F$ is a
non quasiconvex subgroup of $H$. 

Consider now the free product $G = H\ast F$. Since it is well known
that hyperbolic groups are closed under free product, $G$ is
hyperbolic. Moreover, being a free product of torsion-free groups,
it is torsion-free as well.
Let $K = \langle H
\cup aHa\inv\rangle \leq G$ (where $a$ comes from the second factor in
$H\ast F$). It is easy to check that  
\beq
\label{isomfr}
K \cong H\ast H.
\eeq
Indeed, we define a homomorphism $\p:H \ast H \to K$ by sending an
element $h$ from the first factor $H$ into $h \in K$, and an
element $h$ from the second factor $H$ into $aha\inv \in K$. This is
clearly surjective, and injectivity follows from the free product normal form.

We consider now the sequence of embeddings
\beq
\label{seqem}
G = H\ast F \mapright{\theta} H\ast H \mapright{\p} K < H\ast F = G
\eeq
where $\theta$ acts as the identity $H \to H$ with respect to the
first factors and as the inclusion $F \to H$ for the second ones.
Using indices $1$ and $2$ to distinguish generators from different
free factors in the free products, we fix now the finite
generating sets $C$ and $D$ for the $G$ and $H \ast H$, respectively:
$$C = \{ a_1,b_1,t_1,a_2,b_2\},\quad
D = \{ a_1,b_1,t_1,a_2,b_2,t_2 \}.$$
Let $d_C$ and $d_D$ denote the corresponding geodesic metrics on $G$
and $H \ast H$, respectively.

For each $n \geq 0$, let $w_n$ denote the (unique) reduced word over
$\wt{\{ a_2,b_2\}}$ representing the element $t_2^{-n}a_2t_2^n \in F$. It
follows from (\ref{distort2}) and the free product normal form that
\beq
\label{distort3}
d_C(w_n) = 230^n.
\eeq

Now by (\ref{distort1}) we have $d_D(w_n\theta) \leq 2n+1$. We claim that
\beq
\label{distort4}
d_C(w_n\theta\p) \leq 6n+3.
\eeq
Indeed, each generator $a_2,b_2,t_2$ of $H\ast H$ is sent by $\p$ into
$a_2a_1a_2\inv, a_2b_1a_2\inv, a_2t_1a_2\inv$, respectively, and so
$d_D(w_n\theta) \leq 2n+1$ yields $d_C(w_n\theta\p) \leq 3(2n+1) =
6n+3$. Thus (\ref{distort4}) holds. Together with (\ref{distort3}),
this implies that (\ref{seqem}) is not a quasi-isometric embedding.
By the aforementioned theorem of Short, $G$ has a non quasiconvex
subgroup isomorphic to itself. 

\bigskip

If the embedding in such an example can be taken to be uniformly
continuous with respect to some visual metric, we shall have proved
that quasiconvexity cannot be removed from condition (v) in 
Theorem \ref{hpc}(v). But we have no answer yet.

\section{Lipschitz conditions}
\label{lico}

An endomorphism $\p$ of a hyperbolic group $G =
\langle A \rangle$ satisfies an obvious Lipschitz condition
if $G$ is finite. The following two results
provide other instances of Lipschitz conditions.

Given $x \in G$, we denote by $\lambda_x$ the inner
automorphism of $G$ defined by $g\lambda_x = x\inv gx$. 

\bp
\label{hil}
Let $\p$ be an inner automorphism of a 
hyperbolic group $G$ and let $d$ be a visual
metric on $G$. Then $\p$ satisfies a Lipschitz condition.
\ep

\proof
Write $\p = \lambda_x$. Let $d \in V^A(p,\gamma,T)$. 
In view of Proposition \ref{regp}, it suffices to prove that there exists a
constant $Q \in \RR$ such that
\beq
\label{hil1}
(x\inv gx|x\inv hx)_{p}^{A} + Q \geq (g|h)_p^A
\eeq
holds for all $g,h \in G$.

Now
$$\begin{array}{lll}
(x\inv gx|x\inv hx)_{p}^{A}&=&\frac{1}{2}(d_A(p,x\inv gx)+d_A(p,x\inv
hx)-d_A(x\inv gx,x\inv hx))\\
&=&\frac{1}{2}(d_A(xp,gx)+d_A(xp,hx)-d_A(gx,hx))\\
&\geq&\frac{1}{2}(d_A(p,g)-d_A(p,xp)-d_A(g,gx)+d_A(p,h)-d_A(p,xp)\\
&&\hspace{1cm}-d_A(h,hx)
-d_A(g,h)-d_A(g,gx)-d_A(h,hx))\\
&=&(g|h)_p^A-2d(1,x)-d(p,xp)
\end{array}$$
and so (\ref{hil1}) holds for $Q = 2d(1,x)+d(p,xp)$.
\qed

\bl
\label{haq}
Let $\p$ be an automorphism of a 
hyperbolic group $G$ and let $d \in V^A(p,\gamma,T)$ be a visual
metric on $G$. If $\p:(G,d_A) \to
(G,d_A)$ is an isometry, then $\p$ satisfies a Lipschitz condition
with respect to $d$.  
\el

\proof
In view of Proposition \ref{regp}, it suffices to prove that there exists a
constant $Q \in \RR$ such that
\beq
\label{haq1}
(g\p|h\p)_{p}^{A} + Q \geq (g|h)_p^A
\eeq
holds for all $g,h \in G$.

Since
$$\begin{array}{lll}
(g\p|h\p)_{p}^{A}&=&\frac{1}{2}(d_A(p,g\p)+d_A(p,h\p)-d_A(g\p,h\p))\\
&\geq&\frac{1}{2}(d_A(p\p,g\p)+d_A(p\p,h\p)-2d_A(p,p\p)-d_A(g\p,h\p))\\
&=&\frac{1}{2}(d_A(p,g)+d_A(p,h)-d_A(g,h))-d_A(p,p\p) =
(g|h)_p^A-d_A(p,p\p),
\end{array}$$
then (\ref{haq1}) holds for $Q = d_A(p,p\p)$.
\qed

We consider next the particular case of free groups. Given $g \in
F_A$, we denote by $|g|_c$ the {\em cyclic length} of $g$, i.e. the
length of a cyclically reduced conjugate of $g$. 

\bl
\label{fexp}
Let $d \in V^A(p,\gamma,T)$ be a visual
metric on $F_A$ and let $\p$ be an endomorphism of
$G$ satisfying a Lipschitz condition. Then
\beq
\label{fexp1}
|g\p|_c \geq |g|_c
\eeq
for every $g \in F_A$.
\el

\proof
Suppose that there exists some $g \in F_A$ such that $|g\p|_c <
|g|_c$. We may assume that $g$ is cyclically reduced and $g\p = x\inv
hx$ with $h$ cyclically reduced. Hence $|h| < |g|$. For every $n \geq
1$, we have 
$$\begin{array}{lll}
(g^n|g^{n-1})^A_p -
(g^n\p|g^{n-1}\p)^A_p&=&\frac{1}{2}(d_A(p,g^n)+d_A(p,g^{n-1})-d_A(g^{n},g^{n-1})
-d_A(p,g^{n}\p)\\
&&\hspace{1cm}-d_A(p,g^{n-1}\p)+d_A(g^{n}\p,g^{n-1}\p))\\
&\geq&\frac{1}{2}(2d_A(p,g^n)-2d_A(g^n,g^{n-1})
-2d_A(p,g^{n}\p)\\
&=&d_A(p,g^n)-d_A(g,1)-d_A(p,g^{n}\p)\\
&\geq&d_A(1,g^n)-d_A(g,1)
-d_A(1,g^{n}\p)-2d_A(p,1)\\
&=&|g^n|-|g|-|g^n\p|-2|p| = n|g|-|g|-|x\inv h^nx|-2|p|\\
&=&n|g|-|g|-2|x|-n|h|-2|p| \geq n-|g|-2|x|-2|p|.
\end{array}$$ 
Thus there is no
constant $Q \in \RR$ such that
$$(u\p|v\p)_{p}^{A} + Q \geq (u|v)_p^A$$
for all $u,v \in F_A$.
In view of Proposition \ref{regp}, this contradicts $\p$ satisfying a
Lipschitz condition. Therefore (\ref{fexp1}) holds as claimed.
\qed

We can give a complete solution in the case of free group
automorphisms. We denote by $\aut(F_A)$ (respectively $\inn(F_A)$) the
group of automorphisms (respectively inner automorphisms) of
$F_A$. We say that $\p \in \aut(F_A)$ is a {\em permutation
  automorphism} (with respect to $A$) if
$\p|_{\wt{A}}$ is a permutation. Let $\per^A(F_A)$ denote the group 
of all permutation automorphisms of $F_A$ with respect to $A$. 
Since 
\beq
\label{innn}
\lambda_x\p = \p\lambda_{x\p}
\eeq
holds for all $\p \in \aut(F_A)$ and $x \in F_A$, and $\inn(F_A)
\unlhd \aut(F_A)$, $\per^A(F_A) \leq \aut(F_A)$, it follows easily
that
\beq
\label{innn1}
\langle \inn(F_A) \cup
\per^A(F_A) \rangle = \per^A(F_A)\inn(F_A).
\eeq

\bt
\label{fau}
Let $d \in V^A(p,\Gamma,T)$ be a visual
metric on $F_A$ and let $\p \in {\rm Aut}(F_A)$. Then the
following conditions are equivalent:
\bi
\item[(i)] $\p$ satisfies a Lipschitz condition;
\item[(ii)] $\p \in \langle {\rm Inn}(F_A) \cup
{\rm Per}^A(F_A) \rangle = {\rm Per}^A(F_A){\rm Inn}(F_A)$.
\ei
\et

\proof
(i) $\Rw$ (ii). Since $\p$ is an automorphism, it induces a natural
bijection $\oo{\p}$ on the set $C$ of conjugacy classes of $F_A$: if
$[g]$ denotes the conjugacy class of $g$, we set $[g]\oo{\p} = [g\p]$.

Cyclic length extends naturally to conjugacy classes by
$$|[g]|_c = |g|_c.$$
We claim that
\beq
\label{fau1}
|[g]\oo{\p}|_c = |[g]|_c
\eeq
holds for every $g \in F_A$. Indeed, if (\ref{fau1}) fails, the fact
that $\oo{\p}$ is a bijection and there are only finitely many
conjugacy classes of a given cyclic length implies that
$|[g]\oo{\p}|_c < |[g]|_c$ for some $g \in F_A$. Hence
$|g\p|_c < |g|_c$, contradicting Lemma \ref{fexp}. Thus (\ref{fau1})
holds. It follows that
\beq
\label{fau2}
|g\p|_c = |g|_c
\eeq
for every $g \in F_A$.

Let $y \in G$ be such that
$$k = |\{ a \in A : |a\p\lambda_y| = 1\}|$$
is maximum and let $\p' = \p\lambda_y$. 

Suppose that there exists some
$c \in A$ such that $|c\p'| > 1$. We may assume that $A = \{
a_1,\ldots,a_n\}$, $a_i\p' = b_i \in \wt{A}$ for $i = 1,\ldots,k$ and
$c\p' = z\inv d z$ in reduced form. Since $|g\p'|_c = |g\p|_c$ for
every $g \in F_A$, it follows from (\ref{fau2}) that $d \in
\wt{A}$. 

Suppose first that $k > 1$. Then there exist $\varepsilon_1,
\varepsilon_2 \in \{ 1,-1\}$ such that
$(a_1^{\varepsilon_1}ca_2^{\varepsilon_2})\p' = b_1^{\varepsilon_1}z\inv d
zb_2^{\varepsilon_2}$ is a cyclically reduced word. Since
$|a_1^{\varepsilon_1}ca_2^{\varepsilon_2}|_c = 3$ and
$$|(a_1^{\varepsilon_1}ca_2^{\varepsilon_2})\p'|_c = |b_1^{\varepsilon_1}z\inv d
zb_2^{\varepsilon_2}|_c = 3 + 2|z| \geq 5,$$
this contradicts (\ref{fau2}).

Thus $k = 1$. Let $\varepsilon
\in \{ 1,-1\}$ be such that
$(a_1^{\varepsilon}c)\p' = b_1^{\varepsilon}z\inv d
z$ is a reduced word. Write $z = wb_1^m$ in reduced form with $|m|$
maximum. Then
$$(a_1^{\varepsilon}c)\p' = b_1^{\varepsilon}z\inv dz
= b_1^{\varepsilon-m}w\inv dwb_1^m.$$
Similarly to the previous case, $w \neq 1$ implies
$$|(a_1^{\varepsilon}c)\p'|_c = 2|w|+2 > 2 = |a_1^{\varepsilon}c|_c,$$
hence we may assume that $w = 1$. But then $z = b_1^{m}$ and 
$$a_1\p\lambda_{yz\inv} = b_1^{m}b_1b_1^{-m} = b_1,\quad
c\p\lambda_{yz\inv} = c\p\lambda_{y}\lambda_{z\inv} =
c\p'\lambda_{z\inv} = b_1^mz\inv dzb_1^{-m} = d,$$
contradicting the maximality of $k$.

We have thus reached a contradiction in any case following from $k <
|A|$, hence $\p' \in \per^A(F_A)$ and $\p =
\p'\lambda_{y\inv} \in \per^A(F_A)\inn(F_A)$. By (\ref{innn1}),
condition (ii) holds.

(ii) $\Rw$ (i). Permutation automorphisms are clearly isometries of
$(F_A,d_A)$, and satisfying a Lipschitz condition is preserved under
composition. Now the claim 
follows from Proposition \ref{hil} and Lemma
\ref{haq}.
\qed

We note also that Theorem \ref{fau} is far from covering the case of
endomorphisms. For instance, it is easy to check that the monomorphism
$\p$ of $F_A$ defined by $a\p = a^2$ for every $a \in A$ 
satisfies a Lipschitz condition for every visual
metric of the form $d \in V^A(p,\gamma,T)$. 
 
We discuss next what happens when we want to consider arbitrary finite
generating sets. An important role will be played by $\epsilon^A \in
\per^A(F_A)$ defined by $a\epsilon^A = a\inv$ $(a \in A)$.

\bt
\label{fauind}
Let $\p \in {\rm Aut}(F_A)$. Then the
following conditions are equivalent:
\bi
\item[(i)] $\p$ satisfies a Lipschitz condition for every visual
metric on $F_A$;
\item[(ii)] $\p \in {\rm Inn}(F_A)$ or $|A| \leq 1$.
\ei
\et

\proof
(i) $\Rw$ (ii).
In view of Proposition \ref{hil} and Theorem \ref{fau}, it suffices to
show that every $\p \in \per^A(F_A) \setminus \{ 1 \}$ fails condition
(i) when $|A| > 1$.

Suppose first that $a\p \notin \{ a,a\inv\}$ for some $a \in A$. Let
$b = a\p$ and let $A' = (A \setminus \{ b,b\inv\}) \cup \{ u \}$,
where $u = ab$. It is immediate that $A'$ is an alternative basis of
$F_A$. Since $a\p = a\inv u$, we have $|a\p|_c > |a|_c$ in $F_{A'} =
F_A$. Thus $\p$ fails condition (\ref{fau2}) in the proof of Theorem
\ref{fau} when we consider a visual metric $d = \sigma^{A'}_{p,\gamma}$.

Thus we may assume that $a\p = a\inv$ for some $a \in A$. Suppose that
$b\p = b$ for some $b \in A \setminus \{ a \}$. Let $A' = (A \setminus
\{ b\}) \cup \{ u \}$, 
where $u = ab$. Since $u\p = a\inv b = a^{-2}u$, we have  
$|u\p|_c > |a|_c$ in $F_{A'} = F_A$. Thus $\p$ fails condition
(\ref{fau2}) in the proof of Theorem 
\ref{fau} when we consider a visual metric $d =
\sigma^{A'}_{p,\gamma}$.

Finally, we assume that $\p = \epsilon^A$. Fix $b \in
A \setminus \{ a \}$ and let $A' = A \cup \{ u,v\}$, where $u = a^2b$
and $v = a^3b$. Let $w = a^{-2}b\inv a^{-3}b\inv ab$. We claim that
\beq
\label{fauind1}
d_{A'}(1,w^n) = 5n
\eeq
for every $n \in \N$.

Indeed, since $w = a^{-2}v\inv u\inv v$, we have $d_{A'}(1,w) \leq
5$ and so $d_{A'}(1,w^n) \leq 5n$. Assume now that
$$w^n = y_1x_1y'_1x'_1y''_1x''_1y_2x_2y'_2x'_2y''_2x''_2\ldots
y_nx_ny'_nx'_ny''_nx''_nz$$ 
where the $x_i,x'_i,x''_i \in \wt{A'}$ produce the $3n$ occurrences of $b/b\inv$
in the reduced form of $w^n$, and $y_i,y'_i,y''_i,z$ are words on
$\wt{A'}$. Then each $y_i$ must equal $a^{-2}$ after reduction in
$F_A$. Since $x_i$ starts by $b\inv$ and $x''_{i-1}$ ends by $b$ (if
$i > 1$), it follows easily that $y_i$ must have at least two factors
from $\wt{A'}$. Therefore
$$d_{A'}(1,w^n) \geq 3n +2n = 5n$$
and so (\ref{fauind1}) holds.

On the other hand, 
$$w^n\p = (w\p)^n = (a^2ba^3ba\inv b\inv)^n = (uva\inv b\inv)^n$$
yields $d_{A'}(1,w^n\p) \leq 4n$.
Hence
$$(w_n|w_n)_1^{A'} - (w_n\p|w_n\p)_1^{A'} = d_{A'}(1,w^n) -
d_{A'}(1,w^n\p) \geq 5n-4n = n$$
and so there is no
constant $Q \in \RR$ such that
$$(g\p|h\p)_{1}^{A'} + Q \geq (g|h)_1^{A'}$$
for all $g,h \in F_A$.
By Proposition \ref{regp}, $\p$ cannot satisfy a Lipschitz condition
with respect to a 
visual metric $\sigma_{1,\gamma}^{A'}$.

(ii) $\Rw$ (i). If $\p$ is inner, we use Proposition \ref{hil}. If
$|A| \leq 1$, it is easy to see that $\p$ is an isometry of
$(F_A,d_{A'})$ for every finite generating set $A'$ of $F_A$,
therefore we may apply Lemma \ref{haq}.
\qed

The reader may have been surprised by the use of a generating set
which is not a basis in the proof of Theorem \ref{fauind}. The next
two propositions show that we could
have done it for $|A| > 2$, but not for $|A| = 2$:

\bp
\label{basis}
Let $\p \in {\rm Aut}(F_A)$. If $|A| > 2$, then the
following conditions are equivalent:
\bi
\item[(i)] $\p$ satisfies a Lipschitz condition for every visual
metric $d \in V^{A'}(p,\gamma,T)$ on $F_A$, where $A'$ is a basis of $F_A$;
\item[(ii)] $\p \in {\rm Inn}(F_A)$.
\ei
\ep

\proof
By the proof of Theorem \ref{fauind}, it suffices to show that
$\epsilon^A$ satisfies no Lipschitz condition for some visual metric
$\sigma^{A'}_{p,\gamma}$ on $F_A$, 
where $A'$ is a basis of $F_A$. 

Fix distinct $a,b,c \in A$. Once again, let $A' = (A \setminus \{
b \}) \cup \{ u \}$, 
where $u = ab$. It is immediate that $A'$ is an alternative basis of
$F_A$. Since 
$$(uc)\epsilon^A = (abc)\epsilon^A = a\inv b\inv c\inv = a\inv u\inv ac\inv,$$ 
we have $|(uc)\epsilon^A|_c > |uc|_c$ in $F_{A'} =
F_A$. Thus $\epsilon^A$ fails condition (\ref{fau2}) in the proof of Theorem
\ref{fau} when we consider a visual metric $d =
\sigma^{A'}_{p,\gamma}$. Therefore $\epsilon^A$ satisfies no Lipschitz
condition with respect to $\sigma^{A'}_{p,\gamma}$.  
\qed

We fix a canonical basis $A = \{ a,b \}$ for $F_2$ and write
$\epsilon = \epsilon^A$. We note that, similarly to (\ref{innn1}), 
\beq
\label{innn2}
\langle \inn(F_2) \cup
\{ \epsilon \} \rangle = \langle \epsilon \rangle\inn(F_2).
\eeq

Given a basis $\{ u,v \}$ of
$F_2$, we denote by $\mu_{u,v}$ the
automorphism of $F_2$ defined by $a \mapsto u$ and $b \mapsto v$.

\bp
\label{twobasis}
Let $\p \in {\rm Aut}(F_2)$. Then the
following conditions are equivalent:
\bi
\item[(i)] $\p$ satisfies a Lipschitz condition for every visual
metric $d \in V^{A'}(p',\gamma',T)$ on $F_2$, where $A'$ is a basis of $F_2$;
\item[(ii)] $\p \in \langle {\rm Inn}(F_2) \cup
\{ \epsilon \} \rangle = \langle \epsilon \rangle {\rm Inn}(F_2)$. 
\ei
\ep

\proof
(i) $\Rw$ (ii). By the proof of Theorem \ref{fauind}.

(ii) $\Rw$ (i). Let $H = \langle \inn(F_2) \cup
\{ \epsilon \} \rangle$. It is well known that $\aut(F_2)$ is
generated by the Nielsen 
transformations $N = \{ \mu_{b,a}, \mu_{a\inv,b},
\mu_{ab,b} \}$. It is easy to check that
$$\mu_{b,a}\inv\epsilon\mu_{b,a} = \mu_{b,a}\mu_{b\inv,a\inv}
= \epsilon,$$
$$\mu_{a\inv,b}\inv\epsilon\mu_{a\inv,b} = \mu_{a\inv,b}\mu_{a,b\inv}
= \epsilon,$$
$$\mu_{ab,b}\inv\epsilon\mu_{ab,b} = \mu_{ab\inv,b}\mu_{b\inv
a\inv,b\inv} = \mu_{b\inv a\inv b,b\inv} = \epsilon\lambda_b,$$
thus $H \unlhd \aut(F_2)$. 

Let $d' = \sigma^{A'}_{p',\gamma'}$, where $A' = \{u,v\}$ is a basis of
 $F_2$. Then $\mu = \mu_{u,v}$ is an automorphism of $F_2$ such that
 $A\mu = A'$. Thus, for
all $g,h \in F_2$, we have
\beq
\label{twobasis1}
d_{A'}(g,h) = d_{A}(g\mu\inv,h\mu\inv).
\eeq
We claim that
\beq
\label{twobasis2}
(g|h)^{A'}_{p'} = (g\mu\inv|h\mu\inv)^{A}_{p'\mu\inv}
\eeq
holds for
all $g,h \in F_2$.
Indeed, (\ref{twobasis1}) yields
$$\begin{array}{lll}
(g|h)^{A'}_{p'}&=&\frac{1}{2}(d_{A'}(p',g) +
d_{A'}(p',h) - d_{A'}(g,h))\\
&=&\frac{1}{2}(d_{A}(p'\mu\inv,g\mu\inv) +
d_{A}(p'\mu\inv,h\mu\inv) -
d_{A}(g\mu\inv,h\mu\inv)) =
(g\mu\inv|h\mu\inv)^{A}_{p'\mu\inv}. 
\end{array}$$

Now (\ref{twobasis2}) yields
\beq
\label{twobasis3}
(g\epsilon|h\epsilon)^{A'}_{p'} = 
(g\mu\inv\mu\epsilon|h\mu\inv\mu\epsilon)^{A'}_{p'} =
(g\mu\inv\mu\epsilon\mu\inv|h\mu\inv\mu\epsilon\mu\inv)^{A}_{p'\mu\inv}.
\eeq

Since $H \unlhd F_2$, we have $\mu\epsilon\mu\inv \in H$ and so by
Theorem \ref{fau} $\mu\epsilon\mu\inv$ satisfies a Lipschitz condition
with respect to 
any visual metric $d = \sigma^A_{p'\mu\inv,\gamma}$. By
Proposition \ref{regp}, there exists a constant $Q \in \RR$ such that
$$(g\mu\epsilon\mu\inv|h\mu\epsilon\mu\inv)_{p'\mu\inv}^{A} + Q \geq
(g|h)_{p'\mu\inv}^A$$ 
holds for all $g,h \in G$. Together with (\ref{twobasis3}) and
(\ref{twobasis2}), this yields
$$(g\epsilon|h\epsilon)^{A'}_{p'} + Q = 
(g\mu\inv\mu\epsilon\mu\inv|h\mu\inv\mu\epsilon\mu\inv)^{A}_{p'\mu\inv}
+ Q \geq (g\mu\inv|h\mu\inv)^{A}_{p'\mu\inv}
= (g|h)^{A'}_{p'}.$$
By Proposition \ref{regp}, $\epsilon$ satisfies a Lipschitz condition with
respect to $d$. Now the general case follows from Proposition
\ref{hil} and satisfying a Lipschitz condition being preserved by composition.
\qed

We can deduce from the proof of Proposition
\ref{twobasis} a curious result proved by Kapovich, Levitt, Schupp and
Shpilrain in 2007. We recall that $g \in F_A$ is {\em primitive} if it 
belongs to some basis of $F_A$. If $g = a_1a_2\ldots a_n$ in reduced form
$(a_i \in \wt{A})$, the {\em reversal} of $g$ (with respect to $A$) is
defined as $g^R = a_n\ldots a_2a_1$. 

\bc
\label{revcon}
{\rm \cite[Proposition 3.1]{KLSS}}
The conjugacy class of a primitive word of $F_2$ is closed under
reversal.
\ec

\proof
Let $u$ be a primitive word of $F_2$ in reduced form. Then there
exists some $\mu_{u,v} 
\in \aut(F_2)$. Since $H \unlhd F_2$, we have $\mu_{u,v} \epsilon
\mu_{u,v}\inv \in H$. Indeed, it follows from the proof of Proposition
\ref{twobasis} that 
$$\mu_{u,v} \epsilon
\mu_{u,v}\inv = \epsilon\lambda_x$$
for some $x \in F_2$. Hence (\ref{innn}) yields
$$(u^R)\inv = u\epsilon = a\mu_{u,v} \epsilon
 = a\epsilon\lambda_x\mu_{u,v} = a\inv\mu_{u,v}\lambda_{x\mu_{u,v}} =
 u\inv \lambda_{x\mu_{u,v}}$$
and so
$$u^R = u\lambda_{x\mu_{u,v}}.$$
Since a conjugate of a primitive word is itself primitive, we are done.
\qed

In particular, if $u$ is a cyclically reduced primitive word, then
$u^R$, being cyclically reduced, must be a cyclic conjugate of $u$.

Note that Corollary \ref{revcon} does not hold for higher rank (take
e.g. $u = abc$) or for nonprimitive elements of $F_2$ (take e.g. $u =
aba^2b^2$).

\section{Open problems}

The main open problem left by this work
relates to the possibility of
removing quasiconvexity from condition (v) of Theorem \ref{hpc}.

\bq
\label{oq1}
Does every uniformly continuous endomorphism of a hyperbolic group
(with respect to a visual metric) satisfy a H\"older condition? If
not, would it satisfy some other type of condition?
\eq

It would be interesting to extend some of the results in Section
\ref{lico} to hyperbolic, or at least virtually free groups:

\bq
\label{oq3}
When does an automorphism of a hyperbolic
(virtually free) group satisfy a Lipschitz condition?
\eq

\section*{Acknowledgements}

The first author is partially supported by CNPq, PRONEX-Dyn.Syst. and
FAPESB (Brazil). 

\smallskip
 
\noindent
The second author acknowledges support from:
\bi
\item
CNPq (Brazil) through a BJT-A grant (process 313768/2013-7);
\item
the European Regional
Development Fund through the programme COMPETE
and the Portuguese Government through FCT (Funda\c c\~ao para a Ci\^encia e a
Tecnologia) under the project\linebreak
PEst-C/MAT/UI0144/2011.
\ei

\bigskip

{\sc V\'\i tor Ara\'ujo, Universidade Federal da Bahia, Instituto de
  Matem\'atica, Av. Adhemar de Barros, S/N, Ondina, 40170-110
  Salvador-BA, Brazil} 

{\em E-mail address}: vitor.d.araujo@ufba.br

\bigskip

{\sc Pedro V. Silva, Centro de
Matem\'{a}tica, Faculdade de Ci\^{e}ncias, Universidade do
Porto, R. Campo Alegre 687, 4169-007 Porto, Portugal}

{\em E-mail address}: pvsilva@fc.up.pt


\begin{thebibliography}{99}
\bibitem{AS} V.~Ara\'ujo and P.~V.~Silva, Geometric characterizations
  of virtually free groups, preprint, arxiv:1405.5400, 2014.
\bibitem{BH} M.~Bridson and A.~Haefliger, {\em Metric Spaces of
    Non-Positive Curvature}, Grundlehren Math. Wissenschaften, Volume
  319, Springer, New York, 1999. 
\bibitem{Dug} J. Dugundji, {\em Topology}, Allyn and Bacon, 1966.
\bibitem{Fri} A.~H.~Frink, Distance functions and the metrization
  problem, {\em Bull. Amer. Math. Soc.} 43 (1937), 133--142.
\bibitem{GH} E.~Ghys and P.~de~la~Harpe (eds), {\em Sur les Groupes
Hyperboliques d'apr\`es Mikhael Gromov}, Birkhauser, Boston, 1990.
\bibitem{Gro03} M.~Gromov, Random walk in random groups, {\em Geom. Funct.
Anal.} 13(1) (2003), 73--146.
\bibitem{HLV} I.~Holopainen, U.~Lang and A.~V\"ah\"akangas, Dirichlet
  problem at infinity on Gromov hyperbolic metric measure spaces, {\em
    Math. Ann.} 339(1) (2007), 101--134. 
\bibitem{Kap} I.~Kapovich, A non-quasiconvexity embedding theorem for
  hyperbolic groups, {\em Math. Proc. Cambridge Philos. Soc.} 127(3)
  (1995), 461--486.
\bibitem{KB} I.~Kapovich and N.~Benakli, Boundaries of hyperbolic
  groups, In: {\em Combinatorial and geometric group theory}, {\em
    Contemp. Math.} 296, Amer. Math. Soc., Providence, RI, 2002,
  pp. 39--93.  
\bibitem{KLSS} I.~Kapovich, G.~Levitt, P.~Schupp and V.~Shpilrain,
  Translation equivalence in free groups, {\em
    Trans. Amer. Math. Soc.} 359(4) (2007), 1527--1546.
\bibitem{KW} I.~Kapovich and D.~T.~Wise, On the failure of the co-hopf
  property for subgroups of word-hyperbolic groups, {\em Isr.
    J. Math.} 122(1) (2001), 125--147.
\bibitem{LS} R.~C.~Lyndon and P.~E.~Schupp, {\em Combinatorial Group
    Theory}, Springer-Verlag, 1977.
\bibitem{RSel} E.~Rips and Z.~Sela, Structure and rigidity in
  hyperbolic groups I, {\em Geom. Funct. Anal.} 4(3) (1994), 337--371.
\bibitem{Sel} Z.~Sela, Structure and rigidity in (Gromov) hyperbolic
  groups and discrete groups in rank 1 Lie groups II, {\em Geom.
Funct. Anal.} 7(3) (1997), 561-593.
\bibitem{Sho} H.~Short, Quasiconvexity and a theorem of Howson's, In:
  E.~Ghys, A.~Haefliger and A.~Verjovsky (eds.), {\em Group theory
    from a geometrical viewpoint (Trieste, 1990)}, World
  Sci. Publishing, River Edge, NJ, 1991, pp. 168--176.
\bibitem{Sil5} P.~V.~Silva, Fixed points of endomorphisms of virtually
 free groups, {\em Pacific J. Math.} 263(1) (2013), 207--240. 
\bibitem{Vai} J.~V\"ais\"al\"a, Gromov hyperbolic spaces, {\em Expositiones
Math.} 23(3) (2005), 187--231. 
  
\end{thebibliography}
\end{document}